\documentclass{article}
\usepackage[utf8]{inputenc} 
\usepackage[T1]{fontenc}
\usepackage{lipsum}
\usepackage{fullpage}
\title{Notes in the isometric deformations problem in codimension $2$}
\author{Diego N. Guajardo}
\usepackage{graphicx}
\usepackage{amsmath}
\usepackage{amssymb}
\usepackage{pst-node}

\usepackage{tikz-cd} 

\usepackage{bigints}
\usepackage{relsize}
\usepackage{amsthm}
\usepackage{hyperref}
\hypersetup{
    colorlinks=true,
    linkcolor=blue,
    filecolor=magenta,      
    urlcolor=cyan,
    citecolor=blue,
}

\usepackage{tikz-cd}

\makeatletter
\newcommand{\tpitchfork}{%
  \vbox{
    \baselineskip\z@skip
    \lineskip-.52ex
    \lineskiplimit\maxdimen
    \m@th
    \ialign{##\crcr\hidewidth\smash{$-$}\hidewidth\crcr$\pitchfork$\crcr}
  }%
}
\makeatother

\usepackage[sorting=nyt]{biblatex} 

\DeclareFieldFormat[article, inbook]{title}{#1} 
\DeclareFieldFormat[article]{pages}{#1}
\DeclareFieldFormat{journaltitle}{#1\isdot}
\DeclareDelimFormat[bib]{nametitledelim}{: }
\renewbibmacro{in:}{}

\renewbibmacro*{issue+date}{%
  \ifboolexpr{test {\iffieldundef{year}}
              and test {\iffieldundef{issue}}}
    {}
    {\printtext[parens]{%
       \printfield{issue}%
       \setunit*{\addspace}%
       \usebibmacro{date}}}%
  \newunit}

\usepackage[shortlabels]{enumitem}

\newtheorem{thm}{Theorem}[section]

\newtheorem{lema}[thm]{Lemma}
\newtheorem{prop}[thm]{Proposition}
\newtheorem{cor}[thm]{Corollary}

\newtheorem{defn}[thm]{Definition}

\theoremstyle{remark}
\newtheorem{remark}[thm]{Remark}

\usepackage{ifthen} 
\usepackage{xifthen} 

\newcommand{\N}     {\mathbb{N}} 
 
\newcommand{\C}     {\mathbb{C}} 

\newcommand{\R}     {\mathbb{R}}

\renewcommand{\S}{\mathbb{S}}

\newcommand{\de}[2][]{%
    \ifthenelse{\isempty{#1}}
        {\partial_{#2}}
        {\partial_{#2}^#1}
}

\newcommand{\G}[3][]{%
    \ifthenelse{\isempty{#1}}
        {\Gamma_{#2 #3}^#3}
        {\Gamma_{#2 #3}^#1}
}

\newcommand{\D}{\Delta}
\newcommand{\vphi}{\varphi}

\newcommand{\map}[4]{#1:#2^{#3}\rightarrow#4}
\newcommand{\inner}[2]{\langle#1,#2\rangle}

\usepackage{xcolor}

\newcommand{\eref}[1]{(\ref{#1})}

\newcommand{\pref}[1]{Proposition \ref{#1}}

\newcommand{\cref}[1]{Corollary \ref{#1}}

\newcommand{\lref}[1]{Lemma \ref{#1}}

\providecommand{\keywords}[1]
{
  \small	
  \textbf{\textit{Keywords---}} #1
}

\begin{document}

\maketitle\begin{abstract}
     We present a discussion about the local isometric rigidity problem in codimension $2$ with a concrete example. 
     We show the necessity of extending the notions of genuine and honest rigidity in order to have the transitivity property.
     In order to do so, we show the necessity of studying the isometric immersions in semi-Euclidean spaces. 
     We show that this extension comes with a natural type of singularity in the inner product.
\end{abstract}
\hspace{30pt}
\section{Introduction}

Nash's Theorem states that any Riemannian manifold $M^n$ can be isometrically immersed into some Euclidean space. 
The isometric deformation problem is the uniqueness-related question.
Namely, to describe the moduli space of isometric immersions $\map{f}{M}{n}{\R^{n+q}}$ that $M^n$ can have for certain $q$.
In the process, given such an $f$, we need to find a way to meaningfully distinguish it from another one $\map{g}{M}{n}{\R^{n+p}}$.
One way of doing so is using the classical notion of congruence for $p=q$, that is, $f$ and $g$ are 
congruent if there is a rigid motion $T$ of $\R^{n+p}$ such that $f=T\circ g$.
More generally, we say that $f$ and $g$ {\it isometrically extend} if, for some Riemannian manifold $N^{n+r}$ with $0<\ell\leq\min\{p,q\}$, the following diagram is commutative
\begin{equation}\label{diagram isometric extension}
    \begin{tikzcd}[ arrows={-stealth}]
        & & \mathbb{R}^{n+p}\\
        \quad\quad M^n \arrow{urr}{g}\arrow[swap]{drr}{f}     \arrow[hook]{r}{j} %
        & N^{n+r} \arrow[swap]{ur}{G}\arrow{dr}{F} \\
        & & \mathbb{R}^{n+q}
    \end{tikzcd}
\end{equation}
where $j,$ $F$, and $G$ are isometric immersions with $j$ being also an embedding. 
We say that $g$ is a {\it genuine deformation} of $f$ if they do not isometrically extend, not even locally.
This notion was introduced in \cite{DFGenDefSub} and has been extensively used in recent works.

The global isometric deformation problem has been solved for $p+q<\min\{5,n\}$ using a slightly sharper notion than genuine rigidity; some singularities must be allowed in \eref{diagram isometric extension} (see \cite{DGcompactascodim2}, \cite{FGsingular}, and \cite{Sack}).
On the other hand, the local problem has a satisfactory description only for $p=q=1$ and is due to Sbrana \cite{Sbrana} and Cartan \cite{Cartan} in the early $20^{th}$ century.
A modern approach of this result can be found in \cite{DFTinter}.


Recently, there have been some efforts to generalize Sbrana and Cartan's works to codimension $2$.
For example, \cite{DFgenrigcodim2}, \cite{DFThyperEn2}, and \cite{YoSC} present partial results in this matter.
In these works it has been used the following approach: fix a submanifold $\map{f}{M}{n}{\R^{n+q}}$ for $q\in\{1,2\}$ and describe the set of its local {\it genuine deformations} $\map{g}{M}{n}{\R^{n+2}}$.
However, this strategy is non-intrinsic and produces undesired problems.
To solve some of these issues, it was introduced in \cite{FFhyperbolen2} the notion of {\it honest deformation} which is sharper than genuine deformation.
  


As discussed before, some issues have appeared trying to generalize Sbrana and Cartan works.
In addition to that, it is difficult to find non-trivial examples in this theory, for example, it took almost a century to find examples of each family in the Sbrana and Cartan classification; see \cite{DFTinter} and \cite{DFgenrigcodim2}.
This lack of examples makes it difficult to analyze rigidity phenomena in larger codimensions.

In this paper, we analyze the isometric rigidity problem in codimension $2$ using a toy model inspired by \cite{YoSC}. 
The discussion of this example exposes a problem that is not present in the codimension $1$ case. 
Namely, the modern notions of rigidity lack a fundamental property that the classical notion of congruence has; transitivity. 
To be more precise, we found a Riemannian manifold $M^n$ such that any two of its honest deformations $g, \map{\hat{g}}{M}{n}{\R^{n+2}}$ there is a sequence of isometric immersion of $M^n$ in codimension $2$, say $\{g_0,\ldots, g_N\}$ with $g_0=g$ and $g_N=\hat{g}$ such that $\{g_i,g_{i+1}\}$ isometrically extend for $i=0,\ldots, N-1$.
We observe that the ambient space of these immersions may be some semi-Euclidean space $\R^{n+2}_\nu$, that is, $\R^{n+2}$ with a non-degenerate inner product of index $\nu\leq 2$. 

The last discussion suggests that, if we want a transitive concept of rigidity, we may need to analyze the isometric immersions in semi-Euclidean spaces, not only in the Euclidean one.
However, analyzing our toy model in the Lorentz space let us to a natural type of singularity in the inner product of the isometric extension.


\section{Preliminaries}
We start by recalling several concepts used in the isometric rigidity problem. 
Firstly, with the basic notations of linear algebra.

\text{ }

Given a bilinear map $\beta:\mathbb{V}\times \mathbb{U}\rightarrow \mathbb{W}$ between real vector spaces, set
\begin{equation*}
    \mathcal{S}(\beta)=\text{span}\{\beta(X,Y):X\in \mathbb{V}, Y\in \mathbb{U}\}\subseteq \mathbb{W}.
\end{equation*}
The (left) \emph{nullity} of $\beta$ is the vector subspace
\begin{equation*}
    \Delta_\beta=\mathcal{N}(\beta)=\{X\in \mathbb{V}:\beta(X,Y)=0\, ,\,\forall Y\in \mathbb{U}\}\subseteq \mathbb{V}.
\end{equation*}

Assume now that $\mathbb{W}$ has a non-degenerate inner product $\langle\cdot,\cdot\rangle:\mathbb{W}\times \mathbb{W}\rightarrow\mathbb{R}$. 
We say that $\beta$ is $\mathit{flat}$ if
\begin{equation*}
    \langle\beta(X,Y),\beta(Z,W)\rangle=\langle\beta(X,W),\beta(Z,Y)\rangle\quad\forall X,Z\in \mathbb{V}\quad \forall Y,W\in \mathbb{U}.
\end{equation*}


\subsection{Genuine deformations of hypersurfaces in higher codimensions I}\label{seccion de GenDefHyp1}

In this subsection, we comment on two distributions that are central to our work, the nullity and the relative nullity. 
We also recall a notion of isometric rigidity called genuine rigidity.
This concept extends the one used by Sbrana and Cartan for hypersurfaces.

\text{ }

We start by recalling the nullity of a Riemannian manifold $M^n$.

\begin{defn}
    \normalfont The {\it nullity of} $M^n$ {\it at} $x$ is the nullity of the curvature tensor $R$ of $M^n$, that is, the subspace of $T_xM$ given by
\begin{equation*}
    \Gamma(x)=\mathcal{N}(R_x)=\{X\in T_xM: R(X,Y)Z=0,\forall Y,Z\in T_xM\}.
\end{equation*}
The {\it rank of} $M^n$ {\it at} $x$ is defined by $n-\mu$, where $\mu=\dim(\Gamma(x))$. 
\end{defn}
As the results that we are looking for are of local nature and our subspaces are all either kernels or images of smooth tensor fields, we will always  work on each connected component of an open dense subset of $M^n$ where all these dimensions are constant and thus all the subbundles are smooth without further notice. In particular, we assume that $\mu$ is constant and hence the second Bianchi identity implies that $\Gamma$ is a totally geodesic distribution, namely, $\nabla_{\Gamma}\Gamma\subseteq\Gamma$.

For an isometric immersion $f:M^n\rightarrow\mathbb{R}^{n+q}$ we denote by $\alpha^f:TM\times TM\rightarrow T^\perp _fM$ its second fundamental form. 
Let us recall the relative nullity of a submanifold.
\begin{defn}
    \normalfont The {\it relative nullity of} $f$ {\it at} $x$ as the nullity of the second fundamental form $\alpha^f(x)$, that is, the subspace of $T_xM$ given by
    $$\Delta_f(x):=\mathcal{N}(\alpha^f_x).$$
    The {\it rank of} $f$ at $x$ is defined by $n-\nu_f$,  where $\nu_f=\dim(\D_f(x))$.
\end{defn}

Gauss and Codazzi equations show that $\D_f\subseteq\Gamma$ and that $\D_f$ is a totally geodesic distribution.
Moreover, it also shows that $f$ sends the leaves of $\D_f$ into (open subsets of) affine subspaces.
In many circumstances $\D_f=\Gamma$ as shown in \cite{YoChK}, in particular, when $f$ is a nowhere flat hypersurface.

Given two isometric immersions $f:M^n\rightarrow\mathbb{R}^{n+q}$ and $g:M^n\rightarrow\mathbb{R}^{n+p}$, it is useful to work with the vector bundle $W=T^\perp_g M\oplus T^\perp_f M$, in which we define the semi-Riemannian metric with signature $(p,q)$ given by
\begin{equation*}
    \langle(\xi_1,\eta_1),(\xi_2,\eta_2)\rangle=\langle\xi_1,\xi_2\rangle_{T^{\perp}_gM}-\langle\xi_1,\xi_2\rangle_{T^{\perp}_fM}.
\end{equation*}
The bilinear tensor $\beta=(\alpha^g,\alpha^f):TM\times TM\rightarrow W$ is flat with respect to this metric by the Gauss equations of $f$ and $g$.

We say that the pair $\{f,g\}$ {\it extends isometrically} if there exists a Riemannian manifold $N^{n+r}$, an isometric embedding $j:M^n\rightarrow N^{n+r}$ and two isometric immersions $F:N^{n+r}\rightarrow\mathbb{R}^{n+q}$, $G:N^{n+r}\rightarrow\mathbb{R}^{n+q}$ such that $f=F\circ j$ and $g=G\circ j$. That is, the following diagram commutes:
\begin{center}
    \begin{tikzcd}[ arrows={-stealth}]
        & & \mathbb{R}^{n+p}\\
        \quad\quad M^n \arrow{urr}{g}\arrow[swap]{drr}{f}     \arrow[hook]{r}{j} %
        & N^{n+r} \arrow[swap]{ur}{G}\arrow{dr}{F} \\
        & & \mathbb{R}^{n+q}
    \end{tikzcd}
\end{center}
Observe that, in this situation, $\{(G_*\xi,F_*\xi):\xi\in T^{\perp}_jM\}\subseteq \mathcal{S}(\beta)^{\perp}$ is a non-trivial null subbundle of $W$.
Notice that if $\{f,g\}$ and $\{g,h\}$ isometrically extend, then $\{f,h\}$ not necessarily does.
Indeed, see Theorem 9 of \cite{DFTinter} as an example; in that case $\{f_1,F_1\circ f_1\}$ and $\{F_1\circ f_1=F_2\circ f_2,f_2\}$ isometrically extend but $\{f_1,f_2\}$ does not.
This contrasts the notion used by Sbrana and Cartan, which distinguish hypersurfaces by rigid motions of the ambient space.

We say that the pair $\{f,g\}$ is {\it genuine}, or that $g$ is a {\it genuine deformation of $f$} when $f$ is fixed, if there is no open subset $U\subseteq M$ such that $\{f|_U,g|_U\}$ extends isometrically. An isometric immersion $f:M^n\rightarrow\mathbb{R}^{n+q}$ is said to be {\it genuinely rigid} in $\mathbb{R}^{n+p}$ if there is no open subset $U\subseteq M^n$ such that $f|_U$ admits a genuine deformation in $\mathbb{R}^{n+p}$.
If that is not the case, we say that $f$ is {\it genuinely deformable} in $\mathbb{R}^{n+p}$.
In particular, when $f$ is a hypersurface, that $g:M^n\rightarrow\mathbb{R}^{n+p}$ is a genuine deformation of $f$ means that there is no open subset $U\subseteq M^n$ such that $g|_U=h\circ f|_U$, where $h:V\subseteq\mathbb{R}^{n+1}\rightarrow\mathbb{R}^{n+p}$ is some isometric immersion of an open subset $V$ with $f(U)\subseteq V$. 


We denote by $\mathbb{R}^{N}_{\nu}$ the semi-Euclidean space of index $\nu$, that is, $\mathbb{R}^{N}$ with a non-degenerate inner product of index $\nu\leq N$.
All the definitions of this subsection have their natural extensions to the semi-Riemannian context, and we will use them without further mention. 


Finally, let us recall the notion of germs at $x_0\in M^n$. 
Consider $F_{x_0}$ a set of smooth functions defined in a neighborhood of $x_0$ with values in fixed manifold $N^m$.
For $i=1,2,$ let $\map{f_i}{U_i\subseteq M}{n}{N^m}$ be a map of $F_{x_0}$. 
We say that $f_1$ and $f_2$ are identified in the sense of germs at $x_0$ if $f_1|_V=f_2|_V$ for some open subset $V\subseteq U_1\cap U_2$ with $x_0\in V$.
A {\it germ} at $x_0$ is an element of the quotient space defined by this equivalence relation. 
In many aspects, germs behave as maps. 
For example, if $M^n=N^m=\R$ and $F_{0}$ is the set of smooth functions defined in a neighborhood of $0$, then we can sum, multiply, and differentiate germs by simply taking representatives of the germs. 
The composition of germs is also well-defined in many circumstances.
For this reason, it is common to consider germs as maps defined in sufficiently small neighborhood of $x_0$.

\section{On genuine rigidity and transitivity}\label{hyper of rank 3 in cod 2}
In this section, we star with a brief introduction to the problem that we are interested.

\text{ }

As we discussed before, we are interested in the unicity related question associated with Nash's Theorem from a local point of view.
Namely, given a Riemannian manifold $M^n$, a fixed point $x_0\in M^n$, and $p\in\N$, we want to describe the isometric immersions around $x_0$ into flat ambient spaces of codimension $p$.
More precisely, we want to understand the set
$$\mathcal{R}^p_{M^n,x_0}:=\{\map{g}{U\subseteq M}{n}{\R^{n+p}_{\nu}}: g \text{ is an isometric immersion and } U\subseteq M^n \text{ is open with }x_0\in U\}/\sim,$$
where $\sim$ denotes the identification by two relations of equivalence. 
Firstly in the sense of germs at $x_0$, by this means we avoid counting multiple times the same local immersion. 
Moreover, there may be no open subset $U\subseteq M^n$ where all the isometric immersions around $x_0$ are well defined in $U$. 
This shows the necessity of working with germs.
Secondly, we identify maps by rigid motions of the ambient space.
Hence, for simplicity we will always assume that any germ send $x_0$ to $0\in\R^{n+p}$. 
Usually, the manifold or the point will be omitted, so we may write $\mathcal{R}^p_{x_0}$ or $\mathcal{R}^p$ instead.

Recall that, aside of the surfaces, Sbrana and Cartan gave a satisfactory description of $\mathcal{R}^1$ for any Riemannian manifold.
For this reason, we are interested in $\mathcal{R}^2$.

Notice that if $\map{f}{U\subseteq M}{n}{\R^{n+1}_\mu}$ belongs to $\mathcal{R}^1_{M^n,x_0}$ and $h\in\mathcal{R}^1_{\R^{n+1}_\mu,0}$ then the composition $h\circ f$ is in $\mathcal{R}_{M^n,x_0}^2$.
This procedure defines the subset $\mathcal{C}=\mathcal{C}_{x_0}\subseteq\mathcal{R}^2_{x_0}$ of the local compositions around $x_0$.
However, we can argue that $\mathcal{C}$ has not codimension 2 nature nor is intrinsic to $M^n$.
For example, if $f$ has rank at least $3$ then Beez-Killing Theorem shows that $\mathcal{R}^1=\{f\}$, and so $\mathcal{C}$ is completely determined by $\mathcal{R}^1_{\R_{\mu}^{n+1},0}$. 
We also comment that if $f$ has rank at least $4$ then $f$ is genuinely rigid in $\R_{\nu}^{n+2}$ by Theorem 1 of \cite{DTcompii} (the proof can be adapted to semi-Riemannian ambient spaces in codimension $2$), so $\mathcal{R}^2=\mathcal{C}\cong\mathcal{R}^1_{\R_{\mu}^{n+1},0}$ in this case.
\begin{defn}
    \normalfont We say that $g$ is a {\it honest immersion around }$x_0$ if it is not a composition around $x_0$, that is, if it belongs to $$\mathcal{H}^2=\mathcal{H}_{M^n,x_0}^2:=\mathcal{R}^2\setminus\mathcal{C}.$$ 
\end{defn}
\begin{remark}
    In higher codimensions we may need to exclude other compositions in order to have a meaningful concept of honesty.
    Theorem 3.10 of \cite{YoChK} suggests that, in particular, we may need to exclude composition of the form $\map{F}{\hat{M}}{n+q}{\R^{n+p}_\mu}$ and $\map{j}{M}{n}{\hat{M}^{n+q}}$, where $\hat{M}^{n+q}$ is any semi-Riemannian manifold of rank less than $p-q$. 
    Notice that this excludes flat compositions.
\end{remark}

In this work we are interested in analyzing possible structures for $\mathcal{H}^2$.
For this, we will discuss the particular case when $\mathcal{R}^1\neq\emptyset$ and $ \mathcal{H}^2\neq\emptyset$.
As commented before, the rank of $M^n$ must be at most $3$.
To simplify, we assume that $M^n$ has rank $3$ and so $\mathcal{R}^1$ has a unique element $\map{f}{V\subseteq M}{n}{\R^{n+1}}$ (the Lorentzian case is similar). 
We can replace $V$ by $M^n$ since this discussion is from a local perspective. 
In this case, $\mathcal{H}^2$ coincides with the genuine deformations of $f$ (modulo germs and congruence) since $\mathcal{R}^1=\{f\}$.
The next result characterizes the elements of $\mathcal{H}^2$ with an algebraic condition.

\begin{lema}\label{lema g in H iff S(beta)=W}
    Let $\map{f}{M}{n}{\R^{n+1}}$ be a hypersurface of rank $3<n$. 
    Then $g\in\mathcal{H}^2$ if and only if $\mathcal{S}(\beta)=T^{\perp}_gM\oplus T^{\perp}_fM$ where $\beta=(\alpha^g,\alpha^f)$.
\end{lema}
\begin{proof}
    We show first that $\mathcal{S}(\beta)\subseteq T^{\perp}_gM\oplus T^{\perp}_fM$ is a non-degenerate subspace.
    Let $\R^{n+2}_\nu$ be the ambient space of $g$ with $\nu\in\{0,1,2\}$.
    Notice that $\mathcal{S}(\beta)$ is trivially non-degenerate for $\nu=2$.
    Proposition 3.1 of \cite{YoSC} proves the case $\nu=0$, and its proof can be adapted for $\nu=1$.  
    Corollary 17 of \cite{DFGenDefSub} shows the direct part since 
    $$n-3=\dim(\Gamma)\geq\dim(\D_f\cap\D_g)=\dim(\D_\beta)\geq n-\dim(\mathcal{S}(\beta)).$$
    In particular, if $g\in\mathcal{H}^2$ then $\D_g=\Gamma$.
    
    For the converse we simply notice that there are not non-trivial null directions in $\mathcal{S}(\beta)^\perp$.
\end{proof}

It is difficult to characterize $\mathcal{H}^2$ for any rank $3$ Riemannian manifold. 
To illustrate this, let us consider the following example. 
Take $\map{g,h}{L}{2}{\R^3}$ two non-congruent immersions of a surface and consider $\hat{g}:=g\times\text{Id}:L^2\times\R^{n-1}\rightarrow\R^3\times\R^{n-1}=\R^{n+2}$ the product immersion, define $\hat{h}$ in a similar way.
If we intersect $\hat{g}(L^2\times\R^{n-2})$ with a flat hypersurface $F^{n+1}\subseteq\R^{n+2}$ we obtain a hypersurface $\map{f}{M}{n}{F^{n+1}}$. 
We identify $f$ locally with a Euclidean hypersurface. 
Generically, $f$ has rank $3$ and $\hat{h}|_{M^n}$ is a genuine deformation of $f$.
In this case the set $\mathcal{H}^2_{M^n}$ is at least as complicated as $\mathcal{R}^1_{L^2}$.
Observe that $\{\hat{h}|_{M^n}, f\}$ is a genuine deformation, but $\{\hat{h}|_{M^n},\hat{g}|_{M^n}\}$ and $\{\hat{g}|_{M^n},f\}$ are not. 

Since the set $\mathcal{H}^2$ may be to complex for general Riemannian manifolds, we will restrict the nullity to avoid surface-like situations. 
In this simplified setting we will be able to understand isometric rigidity phenomena in higher codimension.
For $T\in\Gamma$ we define the {\it splitting tensor} $C_T:\Gamma^\perp\rightarrow\Gamma^\perp$ as
$$C_T(X):=-(\nabla_XT)_{\Gamma^\perp},$$
where the subindex denotes the orthogonal projection onto $\Gamma^\perp$.
Notice that if $C_T=0$ for all $T\in\Gamma$ then $M^n$ is (locally) a Riemannian product by de Rham's Theorem.
This comment motivates the name of the tensor, it measures how far is $M^n$ of splitting as a product. 
We say that $M^n$ has {\it generic nullity} if there exists $T\in\Gamma$ such that the characteristic polynomial $\psi_{C_T}(z)$ has only simple roots over $\C$, observe that this is an open condition.
We comment that the example of the last paragraph has not generic nullity; the kernel of $C_T$ has rank $2$ or $3$ for all $T\in\Gamma$.

\subsection{Transitivity and genuine deformations}
In this section $\map{f}{M}{n}{\R^{n+1}}$ is a hypersurface of rank $3<n$ with generic nullity and $\mathcal{H}^2\neq\emptyset$.
This defines a family of hypersurfaces big enough to analyze $\mathcal{H}^2$. 

\text{ }

Consider a hypersurface $\map{f}{M}{n}{\R^{n+1}}$ of rank $3$ and generic nullity. 
Let $g\in\mathcal{H}^2=\mathcal{H}^2_{M^n,x_0}$ be an honest immersion.
If the ambient space of $g$ is Euclidean then $f$ and all its genuine deformations in $\R^{n+2}$ are described by Theorem 1.1 of \cite{YoSC}.
We can adapt this result for semi-Euclidean ambient spaces using \lref{lema g in H iff S(beta)=W}.
In particular, we can give a description of $\mathcal{H}^2$ in terms of the maximal parallel flat subbundle of the associated Sbrana bundle.
We present now a summary of how this is done and the results of \cite{YoSC} that will be used in this work.

Let $\map{\pi}{M}{n}{L^3:=M^n/\Gamma}$ be the quotient map of $M^n$ to its leaf space of nullity.
The Gauss map of $f$ descends to the quotient as an immersion $\map{h}{L}{3}{\S^n\subseteq\R^{n+1}}$, and similarly, the expression $\inner{f}{h\circ\pi}$ defines a smooth function $\gamma\in\mathcal{C}(L^3)$ called the {\it support function}.
The Gauss map $h$ and $\gamma$ characterize the hypersurface, that is, $f$ can be parametrized using $h$ and $\gamma$ by means of the Gauss parametrization; see \cite{DGgaussP}.

Since $\D_g=\D_f=\Gamma$, the eigenvectors of the splitting tensor define uniquely (up to order and scaling factors) and intrinsically (independent of $g\in\mathcal{H}^2$) smooth vectors $X_0, X_1, X_2\in\Gamma^{\perp}$ (we are supposing that the eigenvalues of the splitting tensor are real, the complex case is similar).
These vectors diagonalize all the splitting tensors and satisfy $\beta(X_i,X_j)=0$ for $i\neq j$ by Codazzi equation.
Moreover, they descend to $L^3$ as coordinate vectors (after re-scaling factors), that is, there is a chart $(u_0,u_1,u_2)\in\R^3$ such that for $\partial_i:=\partial_{u_i}$ we have
$$\pi_* X_i=\partial_i\circ\pi.$$
This defines a {\it conjugate chart} of the Gauss map $\map{h}{L}{3}{\S^n}$ in the sense that
$$(Q(h))_{ij}=Q_{ij}(h):=\partial^2_{ij}h-\G{i}{j}\partial_jh-\G{j}{i}\partial_ih+g_{ij}h=0,\quad\forall 0\leq i<j\leq 2.$$
Moreover, the support function also satisfies that $Q(\gamma)=0$.

Any $g\in\mathcal{H}^2$ has $3$ distinguished normal vectors
$$\eta_i:=\frac{\alpha^g(X_i,X_i)}{\inner{AX_i}{X_i}}\in T^{\perp}_gM,\quad\forall i,$$
where $A$ is the shape operator associated with the Gauss map of $f$. 
Since $\mathcal{S}(\beta)=T^\perp_gM\oplus T^\perp_fM$ the Gauss equation for $g$ is equivalent to 
\begin{equation}\label{Gauss equation and definition of vphi}
    \inner{\eta_i}{\eta_j}=1+\frac{\delta_{ij}}{\vphi_i},\quad\forall i,j,
\end{equation}
for certain smooth functions $\vphi_i\neq0$, here $\delta_{ij}$ is the Kronecker symbol.
As the normal vectors $\{\eta_i\}_i$ are linearly dependent for dimensional reasons, the matrix $(\inner{\eta_i}{\eta_j})_{ij}$ is singular.
The computation of the determinant shows that
$$\sum_i\vphi_i=-1,$$
and an straightforward computation shows that the linear dependency of the $\eta_i$'s is given by
\begin{equation}\label{eq linear dependency of etas}
    \sum_i\vphi_i\eta_i=0.
\end{equation}
Codazzi equation let us compute the connection $\nabla^{\perp g}$ of the normal bundle $T^{\perp}_gM$ as
\begin{equation}\label{normal connection}
  \left.\begin{array}{rr}
    \nabla^{\perp g}_T\eta_i=0,\quad\,\,\,\,\quad\quad\quad\quad\quad\quad\quad\quad\quad\quad\quad\quad &\forall T\in\Gamma,\,\forall i,\\
    \nabla^{\perp g}_{X_j}\eta_i=\G{j}{i}(\eta_j-\eta_i),\quad\quad\quad\quad\quad\quad\quad\quad\quad &\forall i\neq j,\\
    \nabla^{\perp g}_{X_i}\eta_i=\frac{\partial_i(\vphi_i^{-1})}{2}\vphi_i\eta_i-\vphi_i^{-1}\sum_{j\neq i}\G{i}{j}\vphi_j\eta_j, &\forall i.
  \end{array}\right\}
\end{equation}

Combining \eref{Gauss equation and definition of vphi} with the first equation of \eref{normal connection}, we see that $\vphi$ is constant along the leaves of nullity, so we consider $\vphi_i$ as a smooth map of $L^3$.
Moreover, using the other equations of \eref{Gauss equation and definition of vphi} in a similar way we get that
$$\partial_i\vphi_j=2\G{i}{j}\vphi_j,\quad\forall i\neq j,$$
$$\partial_i\vphi_i=-\partial_i\Big(\sum_{j\neq i}\vphi_j\Big)=-\sum_{j\neq i}2\G{i}{j}\vphi_j,\quad\forall i.$$
These equations naturally define an affine vector bundle $(\R^3\times L^3,\nabla^{\mathcal{S}})$, called the {\it Sbrana bundle}, where the section $\vphi=(\vphi_i)_i$ is parallel.
Furthermore, $g\in\mathcal{H}^2$ is completely determined by $\vphi$, and we explicit this relation by denoting $g=g^{\vphi}$.
For this reason, it is natural to consider the maximal parallel flat subbundle of the Sbrana bundle, say $\mathcal{F}\subseteq L^3\times\R^3$.
Denote by $t(M^n)$, or simply by $t$, the number $(\text{rank}(\mathcal{F})-1)\in\{0,1,2\}$ ($\text{rank}(\mathcal{F})\geq 1$ since $\vphi$ belongs to $\mathcal{F}$). 
If $\overline{x}_0:=\pi(x_0)$ then the map 
\begin{align}\label{bijection Sbrana bundle and deformations}
 \begin{split}
    \mathcal{H}^2_{M^n,x_0}&\longrightarrow O_{\overline{x}_0}\subseteq\mathcal{F}_{\overline{x}_0},\\
    g=g^{\vphi}&\,\,\rightarrow \,\,\vphi(\overline{x}_0),
     \end{split}  
\end{align}
is a natural bijection where 
\begin{equation*}
    O=O_{\overline{x}_0}:=\mathcal{F}_{\overline{x}_0}\cap\{(y_0,y_1,y_2)\in\R^3: y_0+y_1+y_2=-1, y_0y_1y_2\neq 0\}.
\end{equation*}
Notice that $O$ is an open subset of $\mathcal{F}_{\overline{x}_0}\cap\{(y_0,y_1,y_2)\in\R^3: y_0+y_1+y_2=-1\}$ which is diffeomorphic to $\R^{t}$.
We call then the number $t=t(M^n)$ the type of $M^n$, and it measures the size of $\mathcal{H}^2$.

Conversely, suppose that a submanifold $\map{h}{L}{3}{\S^n}$ has a conjugate chart and $\gamma\in\mathcal{C}^\infty(L^3)$ satisfies that $Q(\gamma)=0$. 
Assume also that the maximal parallel flat subbundle $\mathcal{F}$ of its Sbrana bundle satisfies that $O_{\overline{x}_0}\neq\emptyset$ for some $\overline{x}_0\in L^3$.
Then the hypersurface obtained by means of the Gauss parametrization satisfies that $O_{\overline{x}_0}\subseteq\mathcal{H}_{x_0}^2$ for any $x_0$ with $\pi(x_0)=\overline{x}_0$.
Moreover, these sets coincide if the nullity is generic at $x_0$.

\subsection{Transitivity and genuine deformations}
In this section we analyze our toy model. 
This simplified model will present a new phenomenon in the isometric rigidity theory.

\text{ }

Let $\map{f}{M}{n}{\R^{n+1}}$ be a hypersurface of rank $3<n$.
Assume that $M^n$ has generic nullity and fix $g=g^\vphi\in\mathcal{H}^2_{M^n,x_0}$. 
Let us compare $g$ with a distinct honest immersion $\hat{g}=\hat{g}^{\hat{\vphi}}\in\mathcal{H}^2$. 
Denote by $\alpha$ and $\hat{\alpha}$ the corresponding second fundamental forms of $g$ and $\hat{g}$.
Consider 
$$\beta:=(\hat{\alpha},\alpha):TM\times TM\rightarrow T^{\perp}_{\hat{g}}M\oplus T^{\perp}_gM=:W,$$ 
the associated flat bilinear form.
As $\D_g=\Gamma=\D_{\hat{g}}$ and $\alpha(X_i,X_j)=0=\hat{\alpha}(X_i,X_j)$ for $i\neq j$, consider
\begin{equation*}
    \xi_i:=\frac{\beta(X_i,X_i)}{\inner{AX_i}{X_i}}=(\hat{\eta}_i,\eta_i).
\end{equation*} 
By flatness of $\beta$, the set $\{\xi_0,\xi_1,\xi_2\}$ is an orthogonal basis of $\mathcal{S}(\beta)$.
Furthermore, by definition of $\vphi$ and $\hat{\vphi}$ we have that
\begin{equation}\label{norma de xi i}
    \inner{\xi_i}{\xi_i}=\frac{1}{\hat{\vphi}_i}-\frac{1}{\vphi_i},\quad i=0, 1, 2.
\end{equation}
If $\hat{\vphi}_i\neq\vphi_i$ for all $i$, then $g$ is a genuine deformation of $\hat{g}$ since $\mathcal{S}(\beta)^\perp\subseteq W$ is definite.
Hence, assume without loss of generality that $\hat{\vphi}_0=\vphi_0$. 
Consider 
$$\xi:=\hat{\vphi}_1-\vphi_1=\vphi_2-\hat{\vphi}_2\neq 0,$$
which is well defined since $\sum_i\vphi_i=\sum_i\hat{\vphi}_i=-1$.

Consider the line of the Sbrana bundle 
$$L_0:=\{(0,t,-t):t\in\R\}\times L^3\subseteq \R^{3}\times L^3,$$
and observe that $\hat{\vphi}-\vphi=(0,\xi,-\xi)\in L_0$ is a parallel section. 
This implies integrability conditions since
$$2\xi\G{0}{1}=\partial_0(\hat{\vphi}_1-\vphi_1)=\partial_0\xi=\partial_0(\vphi_2-\hat{\vphi}_2)=2\xi\G{0}{2},$$
so 
\begin{equation}\label{(1,0,2) es cero}
    \Gamma_0:=\G{0}{1}=\G{0}{2}.
\end{equation}
Moreover, we also have that
$$\partial_1\xi=\partial_1(\vphi_2-\hat{\vphi}_2)=2\G{1}{2}\xi,$$
and
$$\partial_2\xi=\partial_2(\hat{\vphi}_1-\vphi_1)=2\G{2}{1}\xi$$
which imply the following integrability condition
\begin{equation}\label{(1,2)=(2,1)}
    \partial_1\G{2}{1}=\partial_2\G{1}{2}.
\end{equation}
Straightforward computations show that these are necessary and sufficient conditions for $L_0$ being parallel.
That is, if the $(2,0,1)$-Laplace invariant is zero as in \eref{(1,0,2) es cero}, and the $(2,1)$-Laplace invariant coincides with the $(1,2)$-one as in \eref{(1,2)=(2,1)}, then $L_0$ is parallel. 
In particular, any initial condition $(0,\xi_0,-\xi_0)$ defines a parallel section of the form $(0,\xi,-\xi)$.

Let $I\subseteq\{0,1,2\}$ be the subset of indices $i$ such that $L_i$ is a parallel line of the Sbrana bundle, where $L_1$ and $L_2$ are defined as $L_0$ in a similar way. 
Notice that the size of $|I|$ is related with the type of $M^n$, indeed, any line determined by $I$ is contained in the maximal parallel flat subbundle of the Sbrana bundle.

Given $g\in\mathcal{H}^2$ let $\mathcal{H}_g\subseteq\mathcal{H}^2$ be the subset of the genuine deformations of $g$ in $\mathcal{H}^2$ (or, following \cite{FFhyperbolen2}, the set of honest deformations of $g$ in codimension $2$). 
The next proposition describes the set $\mathcal{H}_g$ in terms of $|I|$.

\begin{prop}\label{proposicion ejemplo de deformaciones en codimension 2}
    Let $\map{f}{M}{n}{\R^{n+1}}$ be a hypersurface of rank $3<n$. 
    Assume that $M^n$ has generic nullity and $g\in\mathcal{H}_{x_0}^2$.
    Then 
    \begin{equation}\label{tres-k}
        \min\{2,|I|\}\leq t(M^n)\leq2.
    \end{equation}
    Furthermore, under the bijection \eref{bijection Sbrana bundle and deformations}, $\mathcal{H}_g$ is identified with $O$ without $|I|$ lines.
    Namely, $\mathcal{H}_g$ is in bijection with:
    \begin{enumerate}[$i)$]
        \item  $O\setminus\{\vphi(\overline{x}_0)\}$ if $|I|=0$;
        \item $O\setminus\{\bigcup_{j\in I}L_j^{g}\}$ if $|I|\neq 0$, where $L_j^{g}\subseteq\mathcal{F}_{\overline{x}_0}^k$ is the affine line $L_j^{g}:=L_j+\vphi(\overline{x}_0)$.
    \end{enumerate}
\end{prop}
\begin{proof}
    Clearly $t\leq 2$ since the Sbrana bundle has rank $3$.
    Moreover, $\vphi$ and the lines defined by $I$ are contained in $\mathcal{F}^k$, this gives the first bound of \eref{tres-k}.

    Take $\tilde{g}\in\mathcal{H}^2$ and let $\tilde{\vphi}$ be the parallel section of the Sbrana bundle that determines $\tilde{g}$.
    We claim that $\{g,\tilde{g}\}$ isometrically extends near $x_0$ if and only if $\vphi_j=\tilde{\vphi}_j$ for a unique $j$, say $j=0$.
    Indeed, we have already proved the direct statement and $(i)$. 
    For the converse, assume that $\vphi_0=\tilde{\vphi}_0$ and thus $0\in I$.
    Let $L=\text{span}\{\eta_0\}\subseteq T^{\perp}_gM$ and $\tilde{L}=\text{span}\{\tilde{\eta}_0\}\subseteq T^{\perp}_{\tilde{g}}M$, where $\xi_0=(\tilde{\eta}_0,\eta_0)$ generates the null space $\mathcal{S}(\beta)\cap\mathcal{S}(\beta)^\perp$.
    Suppose that $\inner{L}{L}\neq 0$ (the case $\inner{L}{L}=0$ will be discussed in the final section of this work).
    Consider the isometry $\tau:L\rightarrow\tilde{L}$ given by $\tau(\eta_0)=\tilde{\eta}_0$ and is trivially parallel since $L$ and $\tilde{L}$ are line bundles. 
    Extend $\tau$ as being the identity in $TM$ and define $\phi_\tau:TM\times(TM\oplus L)\rightarrow \tilde{L}^{\perp}\times L^{\perp}\subseteq T^{\perp}_{\tilde{g}}M\times T^{\perp}_{g}M$ by the formula
    $$\phi_\tau(X,v)=\Big(\big(\tilde{\nabla}_X(\tau(v))\big)_{\tilde{L}^{\perp}},\big(\tilde{\nabla}_Xv\big)_{L^{\perp}}\Big),$$ 
    where $\tilde{\nabla}$ is the ambient space connection and the subindex denotes the orthogonal projection on the respective subspace. 
    We want to use Proposition 12 of \cite{FGsingular} to prove that $\{g,\tilde{g}\}$ isometrically extends around $x_0$.
    
    Notice that $\phi_\tau(\Gamma,TM)=\phi_\tau(X_0,TM)=0$ since $\beta(X_0,X_i)=0$ for $i\neq 0$ and $\beta(X_0,X_0)\in\tilde{L}\times L$.
    Moreover, by \eref{eq linear dependency of etas} and \eref{normal connection} we have that $\phi_\tau(\Gamma,L)=0$ and 
    \begin{equation}\label{B}
        \nabla_{X_0}^{\perp g}\eta_0=\frac{\de{0}(\vphi_0^{-1})}{2}\vphi_0\eta_0-\Gamma_0\vphi_0^{-1}(\vphi_1\eta_1+\vphi_2\eta_2)=\Big(\frac{\de{0}(\vphi_0^{-1})}{2}\vphi_0+\Gamma_0\Big)\eta_0\in L.
    \end{equation}
    By an analogous formula for $\tilde{\eta}_0$ we conclude that $\phi_\tau(X_0,L)=0$. 
    Hence, $\phi_\tau(\text{span}\{X_0\}+\Gamma,TM\oplus L)=0$.
    Using \eref{normal connection} we notice that
    $$\phi(X_i,\eta_0)=\G{i}{0}\big((\tilde{\eta}_i)_{\tilde{L}^{\perp}},(\eta_i)_{L^{\perp}}\big)=\G{i}{0}\inner{AX_i}{X_i}^{-1}\phi(X_i,X_i),\quad \forall i\neq0.$$
    Define then
    \begin{equation}\label{lambda}
        \lambda:=\eta_0-\frac{\G{1}{0}}{\inner{AX_1}{X_1}}X_1-\frac{\G{2}{0}}{\inner{AX_2}{X_2}}X_2\in TM\oplus L,
    \end{equation}
    and consider the line bundle $\Lambda=\text{span}(\lambda)\subseteq TM\oplus L$. 
    The computations above show that $\phi(TM,\Lambda)=0$.
    Let $\R^{n+2}_\nu$ be the ambient space of $g$ and define $G:N^{n+1}\subseteq\Lambda\rightarrow\R^{n+2}_\nu$ given by
    \begin{equation}\label{G}
        G(t\eta_p)=g(p)+t\eta_p,
    \end{equation}
    where $N^{n+1}$ is a sufficiently small neighborhood of the zero section $j:M^n\rightarrow\Lambda$ in order to $G$ being an immersion, and define $\Tilde{G}$ in a similar way.
    Proposition 12 of \cite{FGsingular} shows that $G$ and $\Tilde{G}$ induce the same metric on $N^{n+1}$, so $\{g,\tilde{g}\}$ isometrically extends since $g=G\circ j$ and $\Tilde{g}=\Tilde{G}\circ j$. 
    This proves $(ii)$. 
    We observe that the inner product of $N^{n+1}$ is Riemannian when the inner product of $L$ is positive definite (that is, if $1+\vphi_0^{-1}>0$), and Lorentzian on the contrary (when $1+\vphi_0^{-1}<0$).
    \end{proof}
    \begin{remark}
    In the last proof the pair $\{g,\tilde{g}\}$ isometrically extends as a Sbrana-Cartan hypersurface of the continuous or discrete families; see in \cite{DFTinter} the definition of these families for the Riemannian case, for semi-Riemannian metrics the definitions are analogous.
\end{remark}

    \begin{cor}\label{cor asd}
        Under the hypothesis of the last result, if $|I|\geq 2$ and $\hat{g}\in\mathcal{H}_g$, then we have that:
        \begin{enumerate}[$i)$]
            \item There exists $g_1\in\mathcal{H}^2$ such that $\{g,g_1\}$ and $\{g_1,\hat{g}\}$ isometrically extend or;
            \item There are $g_1,g_2\in\mathcal{H}^2$ such that $\{g,g_1\}$, $\{g_1,g_2\}$, and $\{g_2,\hat{g}\}$ isometrically extend.
        \end{enumerate}
    \end{cor}
    \begin{proof}
    Let $\hat{\vphi}$ be the parallel section of the Sbrana bundle that defines $\hat{g}\in\mathcal{H}_g$. 
    Assume without loss of generality that $\{0,1\}\subseteq I$.
    Denote by $\ell_0$ and $\ell_1$ the parallel sections of $L_0$ and $L_1$, respectively, with $\ell_0(x_0)=(0,1,-1)$ and $\ell_1(x_0)=(-1,0,1)$.
    Let $a,b\in\R$ be the unique constants such that 
    $$\vphi-\hat{\vphi}=a\ell_0-b\ell_1$$
    Consider the parallel section $\vphi^1=\vphi+b\ell_1=\hat{\vphi}+a\ell_0$.
    If his entries are non-zero at $x_0$ then it defines an element of $g_1\in\mathcal{H}^2$ which satisfies $(i)$ since $\vphi^1_0=\vphi_0$ and $\vphi^1_1=\hat{\vphi}_1$, the same argument works for the parallel section $\vphi-a\ell_0=\hat{\vphi}-b\ell_1$.
    Hence, it remains to analyze the case in which both parallel sections have a zero entry at $x_0$, a straightforward computation shows that in this case 
    $\vphi-\hat{\vphi}=\vphi_2(x_0)(\ell_1-\ell_0)$.
    Define, for example, $\vphi^1:=\vphi+\frac{\vphi_2(x_0)}{2}\ell_0$ and $\vphi^2=\vphi^1-\vphi_2(x_0)\ell_1$.
    These parallel sections define $g_1,g_2\in\mathcal{H}^2$ satisfying $(ii)$.
\end{proof}
\begin{remark}
    Even when the ambient space of $g$ and $\hat{g}$ is the Euclidean ambient space $\R^{n+2}_0=\R^{n+2}$, this is not necessarily the case for $g_1$.
\end{remark}

There are many hypersurfaces of rank $3$ with $|I|\neq 0$. 
For example, take a Sbrana-Cartan hypersurface $\map{F}{\hat{M}}{n+1}{\R_{\mu}^{n+2}}$ of the continuous family, and let $\map{G_t}{\hat{M}}{n+1}{\R_{\nu_t}^{n+2}}$ be the genuine deformations of $F$ parametrized by an open subset $t\in J\subseteq\R$. 
Intersect $F_2:=F$ with a flat hypersurface $\map{F_1}{U\subseteq\R}{n+1}{\R_{\mu}^{n+2}}$ that is not contained on a hyperplane.
The following diagram describes our situation
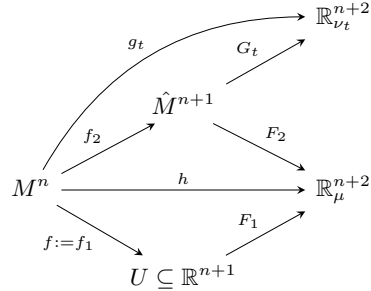
\begin{equation*}
    \begin{tikzcd}[ arrows={-stealth}]
        & & \mathbb{R}^{n+2}_{\nu_t}\\    
        &  \hat{M}^{n+1}\arrow{dr}{F_2}\arrow{ur}{G_t}  \\
        M^n\arrow[bend left=30]{uurr}{g_t} \arrow{ur}{f_2}\arrow[swap]{dr}{f:=f_1}\arrow{rr}{h}      
        & & \mathbb{R}_{\mu}^{n+2}  \\
        & U\subseteq\mathbb{R}^{n+1} \arrow{ur}{F_1}
    \end{tikzcd}
\end{equation*}
Let use the notations of the above diagram. 
Generically, $f$ is a rank $3$ hypersurface of generic nullity and $g_t$ is a genuine deformation of $f$. 
Observe that $\{g_{t_1},g_{t_2}\}$ isometrically extends if $t_1\neq t_2$ since $g_t:=G_t\circ f_2$ for $t\in J$.
This shows that $f$ has $|I|\neq0$. 

Notice that the set $\mathcal{H}_g$ depends of the immersion $g$, in particular, it is not intrinsic.
However, the last results suggests a natural equivalence relation $\sim$. 
Given $g, \hat{g}\in\mathcal{R}^2_{x_0}$, we say that they are related if there exists a sequence of isometric immersions $g_0,\ldots ,g_N$ such that $g_0=g$, $g_N=\hat{g}$, and $\{g_{i-1},g_{i}\}$ isometrically extends around $x_0$ for $i=1,\ldots,N$.
\cref{cor asd} shows the following result.
\begin{thm}
    In the situation of \pref{proposicion ejemplo de deformaciones en codimension 2}, $\mathcal{H}^2/\sim$ is in a natural bijection to an open subset of $\R^{t-\min\{|I|,2\}}$, where $t=t(M^n)$ is the type of $M^n$.
\end{thm}


\section{Case \texorpdfstring{$\inner{L}{L}=0$}{TEXT}}
In the proof of \pref{proposicion ejemplo de deformaciones en codimension 2} we skip the case in which $\inner{L}{L}=0$.
Analyzing this case will help to understand a phenomena that appears by allowing semi-Euclidean ambient spaces. 
The description of this case suggests that we may need to allow a type of singularity in the inner product. 

\text{ }

Recall the situation of \pref{proposicion ejemplo de deformaciones en codimension 2}, we have two honest immersions $g,\tilde{g}\in\mathcal{H}^2$ with $\vphi_0=\tilde{\vphi}_0$ where $\vphi$ and $\tilde{\vphi}$ are the parallel sections of the Sbrana bundle that define $g$ and $\tilde{g}$, respectively.
Then, the normal subbundle $L=\text{span}\{\eta_0\}$ and $\Tilde{L}=\text{span}\{\tilde{\eta}_0\}$ are isometric and parallel. 
The inner product in $L$ (and $\Tilde{L}$) is determined by 
$$\inner{\eta_0}{\eta_0}=\inner{\Tilde{\eta}_0}{\tilde{\eta}_0}=1+\frac{1}{\vphi_0}.$$
We proved that $g$ and $\tilde{g}$ isometrically extend when the inner product non-degenerate on $L$, that is, when $\vphi_0\neq-1$. 
To be more precise, we explicit its isometric extensions; the maps $G$ and $\Tilde{G}$ defined by \eref{G} (where $\Tilde{G}$ is defined analogously).
It remains to analyze the case $\vphi_0=-1$ ($\vphi_0$ is not necessarily identically $-1$).
Observe that this case only happens when the ambient space of $g$ and $\tilde{g}$ is the Lorentz space $\R^{n+2}_1$.

The maps $G$ and $\Tilde{G}$ are well defined even in the case that $\vphi_0=-1$.
Hence, to analyze the case $\vphi_0=-1$, it is natural to study the immersions $G$ and $\tilde{G}$ in this situation. 
However, we should expect some type of the singularity on the metric tensor induced by these immersions. 
Indeed, around points where $\inner{L}{L}=0$,
the inner product of $N^{n+1}$ induced by $G$ (or $\tilde{G}$) may change its signature.
We also observe that we may have $\vphi_0\equiv-1$ in a neighborhood of $x_0$, hence, we cannot exclude this case by an argument of density.

\begin{lema}
    Both immersions $G$ and $\tilde{G}$ are isometric in the sense that $T:=G^*h=\Tilde{G}^*h$, where $h$ is the canonical inner product of the Lorentz space. 
    Moreover, $T$ is a Riemannian or Lorentzian inner product along $\Lambda$ where $(1+\vphi_0^{-1})$ is positive or negative, respectively.
    However, $T$ degenerate along $\Lambda$ where $\vphi_0=-1$.
\end{lema} 
\begin{proof}
    We have already proved this for $\vphi_0\neq -1$.
    Assume that $\vphi_0(\overline{x}_0)=-1$ where $\overline{x}_0=\pi(x_0)\in L^{3}$ is fixed.
    Consider an smooth function $\theta:M^n\rightarrow\R$ and the section $\xi=\theta\lambda$, where $\lambda$ is given by \eref{lambda}. 
    It is enough to show that 
    \begin{equation}\label{A}
        A:=h((G\circ\xi)_*X,(G\circ\xi)_*Y)-h((\tilde{G}\circ\xi)_*X,(\tilde{G}\circ\xi)_*Y)=0.
    \end{equation}
    Write $\lambda=\eta_0+Y$, where $Y=-\frac{\G{1}{0}}{\inner{AX_1}{X_1}}X_1-\frac{\G{2}{0}}{\inner{AX_2}{X_2}}X_2$, and notice that 
    \begin{align*}
        (G\circ\xi)_*X&=X+\tilde{\nabla}_X\xi=X+(\tilde{\nabla}_X\xi)_{TM}+(\tilde{\nabla}_X\xi)_{T^{\perp}_gM}=\\
        &=(X-\theta A_{\eta_0}X+d\theta(X)Y+\theta\nabla_XY)+(d\theta(X)\eta_0+\theta\nabla_X^{\perp}\eta_0+\theta\alpha^g(X,Y))=\\
        &=T_g(X)+N_g(X)\in TM\oplus T^{\perp}_gM,
    \end{align*}
    where $A_{\eta_0}$ is the shape operator associated to $\eta_0\in T_{g}^{\perp}M$. Notice that $T_g(X)=T_{\tilde{g}}(X)$ since $A_{\eta_0}=A_{\tilde{\eta}_0}$. 
    Hence, to verify \eref{A}, it is enough to show that 
    $$A=h(N_g(X),N_g(Y))-h(N_{\Tilde{g}}(X),N_{\Tilde{g}}(Y))=0,$$
    using \eref{normal connection} and \eref{B} we verify that $N_g(X)=a(X)\eta_0$ for some $1$-form $a$ which only depends on $\vphi_0=\Tilde{\vphi}_0$, and similarly $N_{\Tilde{g}}(X)=a(X)\Tilde{\eta}_0$.
    This concludes the proof since $\eta_0$ and $\Tilde{\eta}_0$ are isometric.
\end{proof}

\printbibliography

IMPA – Estrada Dona Castorina, 110

22460-320, Rio de Janeiro, Brazil

{\it E-mail address}: {\tt diego.navarro.g@ug.uchile.cl}
\end{document}